\documentclass[10pt,leqno]{amsart}%
\usepackage{amsmath}
\usepackage{amsfonts}
\usepackage{amssymb}
\usepackage{graphicx}%
\usepackage{color}
\usepackage{hyperref}

\providecommand{\U}[1]{\protect\rule{.1in}{.1in}}
\newtheorem{theorem}{Theorem}

\newtheorem{corollary}[theorem]{Corollary}

\newtheorem{lemma}[theorem]{Lemma}

\newtheorem{remark}[theorem]{Remark}

\newtheorem{thm}{Theorem}
\newtheorem*{acknowledgement*}{Acknowledgement}

\newcommand{\tr}{\operatorname{tr}}

\newcommand{\kp}[1]{k_p\left(#1\right)}
\newcommand{\Kp}[1]{\mathcal{K}_p\left(#1\right)}

\newcommand{\tnn}{\tilde N_{\times/{}}}

\newcommand{\Hess}{\operatorname{Hess}}

\newcommand{\sect}{\operatorname{Sect}}

\newcommand{\nsect}{{}^N \operatorname{Sect}}

\newcommand{\vol}{\operatorname{Vol}}

\newcommand{\tu}{\tilde{u}}
\newcommand{\tir}{\tilde{r}}
\newcommand{\tq}{\tilde{q}}
\newcommand{\tv}{\tilde{v}}

\newcommand{\tn}{\tilde{N}}
\newcommand{\tm}{\tilde{M}}

\newcommand{\dive}{\operatorname{div}}

\newcommand{\mmetr}[2]{\left\langle #1,#2\right\rangle_M}
\newcommand{\nmetr}[2]{\left\langle #1,#2\right\rangle_N}
\newcommand{\hsmetr}[2]{\left\langle #1,#2\right\rangle_{HS}}

\newcommand{\nr}[2]{{}^N r(#1,#2)}
\newcommand{\dist}{\operatorname{dist}}

\begin{document}

\title[Comparison theorem for $p$-harmonic maps]{A general comparison theorem for $p$-harmonic maps in homotopy class}
\date{\today}
%

\author{Giona Veronelli}
\address{D\'epartement de Math\'ematiques\\
Universit\'e de Cergy-Pontoise\\
Site de Saint Martin
2, avenue Adolphe Chauvin\\
95302 Cergy-Pontoise Cedex,
France.}
\email{giona.veronelli@gmail.com}

\subjclass[2010]{58E20, 53C43}

\keywords{$p$-harmonic maps, homotopy classes of maps, $p$-parabolicity}

\begin{abstract}
We prove a general comparison result for homotopic finite $p$-energy $C^{1}$ $p$-harmonic maps $u,v:M\to N$ between Riemannian manifolds, assuming that $M$ is $p$-parabolic and $N$ is complete and non-positively curved. In particular, we construct a homotopy through constant $p$-energy maps, which turn out to be $p$-harmonic when $N$ is compact. Moreover, we obtain uniqueness in the case of negatively curved $N$. This generalizes a well known result in the harmonic setting due to R. Schoen and S.T. Yau.
\end{abstract}
\maketitle

\section{Introduction and main result}

In the mid 1960's, J. Eells and J.H. Sampson extended the notion of harmonicity from real-valued functions to manifold-valued maps, \cite{ES}. The topological relevance of harmonic maps, already visible in the seminal paper \cite{ES}, became clear in the works by P. Hartman, L. Lemaire, R. Hamilton and others authors. See e.g. \cite{Ha}, \cite{EL}, \cite{Ham}. Notably, R. Schoen and S.T. Yau developed the theory and topological consequences of harmonic maps with finite energy from a non-compact domain, \cite{SY-Helv}, \cite{SY-Topo}, \cite{SY-LN}.\\
A natural extension of the concept of a harmonic map is that of a $p$-harmonic map. To this end, a great deal of work has been done by B. White, \cite{Wh}, R. Hardt and F.-H. Lin, \cite{HL}, and S.W. Wei, \cite{Wei-indiana}, \cite{Wei-conference}. In particular, it is known that $p$-harmonic maps give information on the higher homotopy 
groups, and on the homotopy class of higher energy maps between Riemannian manifolds, \cite{Wei-indiana}, \cite{PV}. In view of these topological links we are led to understand which of the well known results holding in the harmonic case can be extended to the non-linear setting. \\
First, we recall that a $C^1$ map $u:(M,\mmetr{}{})\to(N,\nmetr{}{})$ between Riemannian manifolds is said to be $p$-harmonic, $p>1$, if its $p$-tension 
field $\tau_p u$ vanishes everywhere, i.e. if $u$ satisfies the non-linear system
\begin{equation}\label{Itau_p}
\tau_{p}u=\operatorname{div}\left(  \left\vert du\right\vert ^{p-2}%
du\right)=0.
\end{equation}
Here, $du\in T^{\ast}M\otimes u^{-1}TN$ denotes the differential of $u$ and the
bundle $T^{\ast}M\otimes u^{-1}TN$ is endowed with its Hilbert-Schmidt scalar
product $\hsmetr{}{}$. Moreover, $-\operatorname{div}$ stands for the formal adjoint of the
exterior differential $d$ with respect to the standard $L^{2}$ inner product
on vector-valued $1$-forms. Observe that, when $N=\mathbb R$, $\tau_p$ coincides with the standard
 $p$-laplace operator $\Delta_p$. Clearly, in general equality (\ref{Itau_p}) has to be considered in the weak sense, i.e.
\[
\int_M\hsmetr{|du|^{p-2}du}{\eta}=0
\] 
for every smooth compactly supported $\eta\in \Gamma(T^{\ast}M\otimes u^{-1}TN)$. In case $p=2$, the non-linear factor $|du|^{p-2}$ disappears and the $2$-harmonic map is simply
called harmonic. Furthermore, by standard elliptic regularity, a harmonic map $u$ is necessarily smooth.\\
When the target manifold $N$ is compact with nonpositive sectional curvatures, the existence of a $p$-harmonic representative in the homotopy class of a given map has been proved by S.W. Wei, \cite{Wei-conference}. Namely, elaborating on ideas of F. Burstall, \cite{Burstall-London}, he gave the following result, a detailed proof of which is contained in the next section.

\begin{theorem}[Wei]\label{th_Wei-existence}
Let $M$ be a complete Riemannian $m$-dimensional manifold and $N$ a compact manifold with $\nsect\leq 0$. Then any continuous (or more generally $W^{1,p}$) map $f:M\to N$ of finite $p$-energy, $1<p<\infty$, can be deformed to a $C^{1,\alpha}$ $p$-harmonic map $u$ minimizing $p$-energy in the homotopy class.
\end{theorem}

According to Theorem \ref{th_Wei-existence}, $p$-harmonic maps can be considered as ``canonical'' representatives of homotopy class of maps with finite $p$-energy. Hence, one is naturally led to investigate such a space, in particular inquiring how many $p$-harmonic 
representatives can be found 
in a given homotopy class.\\
A first uniqueness result in this direction was obtained by Wei, \cite{Wei-indiana}, for smooth $p$-harmonic maps defined on 
compact $M$, generalizing a previous result for $p=2$ due to P. Hartman, \cite{Ha}. An interesting task is then to detect similar results for complete non-compact 
manifolds. In the harmonic setting, the most important result is represented by the following celebrated theorem due to Schoen and Yau, \cite{SY-Topo}, \cite{SY-LN}.
 
\begin{theorem}[Schoen-Yau]\label{Ith_SY_Topo}
Let $M$ and $N$ be complete Riemannian manifolds with $\vol M<\infty$.
\begin{itemize}
 \item[i)] Let $u:M\to N$ be a harmonic map of finite energy. If $\nsect <0$, there's no other harmonic map of finite energy homotopic to $u$ unless $u(M)$ 
is contained in a geodesic of $N$.
\item[ii)] If $\nsect\leq 0$ and $u,v:M\to N$ are homotopic harmonic maps of 
finite energy, then there is a smooth one-parameter family $u_t:M\to N$, of harmonic maps with $u_0=u$ and $u_1=v$. Moreover, for each 
$q\in M$, the curve $\{u_t(q):t\in\mathbb R\}$ is a constant (independent of $q$) speed parametrization of a geodesic. Finally, the map $M\times\mathbb R\to N$ given by $(q,t)\mapsto u_t(q)$ is harmonic with respect to the product metric on $M\times\mathbb R$.
\end{itemize}
\end{theorem}  

Theorem \ref{Ith_SY_Topo} has been subsequently improved by S. Pigola, M. Rigoli and A.G. Setti replacing the finite volume condition $\vol M<+\infty$ with the parabolicity of $M$, \cite{PRS-MathZ}. We recall that a manifold $M$ is said to be $p$-parabolic if any bounded above weak solution $\varphi$ of $\Delta_p\varphi\geq 0$ is necessarily constant. It is well known that this is just one of the several equivalent definitions of $p$-parabolicity; see \cite{Holo-Dissertation}, \cite{Troy}, \cite{PST} and Theorem \ref{weak_KNR} below. Moreover, a $2$-parabolic manifold is called simply parabolic.\\
As a matter of fact, an inspection of Schoen and Yau's proof shows that they strongly use the property of ($2$-)harmonic maps to be 
($2$-)subharmonic, once composed with a convex function. It turns out that, in general, this is not true if $p\neq 2$, \cite{V}. Hence, one is led to follow different paths in order to to deal with the non-linear analogous of 
Schoen and Yau's result. In this direction, some progresses in the special
situation of a single map homotopic to a constant has been made in \cite{PRS-MathZ}, where the authors managed to overcome the bad behaviour of the composition of $p$-harmonic maps with convex functions introducing special composed vector fields and applying to this latter a global form of
the divergence theorem in non-compact settings, due to V. Gol'dshtein and M. Troyanov, which goes under the name of Kelvin-Nevanlinna-Royden criterion, \cite{GT}; see Theorem \ref{weak_KNR} below. The first attempt to consider two non-constant maps has been made in \cite{HPV}. In that paper, the case
$N=\mathbb{R}^{n}$ has been considered. According to \cite{PRS-MathZ}, if $M$ is
$p$-parabolic, then every $p$-harmonic map $u:M\rightarrow\mathbb{R}^n$ with
finite $p$-energy $\left\vert du\right\vert \in L^{p}\left(  M\right)  $ must
be constant. However, using the very special structure of $\mathbb{R}^{n}$ the authors were able to extend this conclusion obtaining that two finite $p$-energy maps $u,v:M\rightarrow\mathbb{R}^n$ differs by a constant provided $\tau_p u=\tau_p v$. Since $\mathbb R^n$ is contractible, in this situation all maps are trivially homotopic. Consequently, this result can be seen as a special case in the comprehension of general comparison theorems for homotopic $p$-harmonic maps in the spirit of Theorem \ref{Ith_SY_Topo}. However, though the procedure is non trivial due to the non linearity of $\tau_p$, thanks to the good properties of $\mathbb R^n$, i.e. in particular the standard way to compare $p$-tension fields at different points through their difference, the problem in \cite{HPV} was somehow reduced to that of a single map. Nevertheless, it is worth noting that in this case the Kelvin-Nevanlinna-Royden criterion has been used in his full power in the sense that the the vector field $X$ introduced has a non-trivial negative part. This leads to employ a new limit procedure which turns out to be useful in our investigation; see also Theorem 15 in \cite{VV}.\\

The main achievement of this paper is the following very general comparison result for homotopic $p$-harmonic maps which represents an analogue of Theorem \ref{Ith_SY_Topo} for $p\geq 2$.  

\begin{thm}\label{th_SYp}
Suppose $M$ is $p$-parabolic, $p\geq2$, and $N$ is complete. 
\begin{itemize}
\item[i)] If $\nsect<0$ and $u:M\to N$ is a $C^{1}$ $p$-harmonic map of finite $p$-energy, then there's no other $p$-harmonic map of finite $p$-energy homotopic to $u$ unless $u(M)$ is contained in a geodesic of $N$. 
\item[ii)]
If $\nsect\leq 0$ and $u,v:M\to N$ are homotopic $C^{1}$ $p$-harmonic maps of finite $p$-energy, then there is a continuous one-parameter family of maps $u_t:M\to N$ with $u_0=u$ and $u_1=v$ such that the $p$-energy of $u_t$ is constant (independent of $t$) and for each $q\in M$ the curve $t\mapsto u_t(q)$, $t\in[0,1]$, is a constant (independent of $q$) speed parametrization of a geodesic. Moreover, if $N$ is compact, $u_t$ is a $p$-harmonic maps for each $t\in[0,1]$. 
\end{itemize}
\end{thm}

\begin{acknowledgement*}
This work forms part of the author's doctoral thesis. The author is deeply grateful to his advisor Stefano Pigola and to Alberto Setti for their constant encouragement and their useful advice.
\end{acknowledgement*}

%

\section{Proofs of the results}

Since the arguments in \cite{Wei-conference} are only sketched, and since the proof of this result plays a key role in the development of the paper, we provide a detailed proof of Theorem \ref{th_Wei-existence}. To this end, we adapt to the case $p\geq 2$ the proof given by Burstall for $p=2$, \cite{Burstall-London}. 

\begin{proof}[Proof (of Theorem \ref{th_Wei-existence})]
For the ease of notation, throughout all the proof we will keep the same set of indeces each time we will extract a subsequence from a given sequence.\\ 
Consider an exhaustion $\{M_k\}_{k=1}^{\infty}$ of $M$, i.e. a sequence such that, for each $k$, $M_k$ is a compact manifold with boundary, $M_k\subset\subset M_{k+1}$ and $\cup_{k=1}^{\infty}M_k=M$. Define $\mathcal{H}_f$ as the space of $W_{loc}^{1,p}(M,N)$ maps $v$ such that $v|_{M_k}$ and $f|_{M_k}$ have the same $1$-homotopy type, i.e.
\[
\mathcal H_f := \{v\in W_{loc}^{1,p}(M,N) : \forall k\geq 1,\ \left(v|_{M_k}\right)_{\sharp}\textrm{ is conjugated to }\left(f|_{M_k}\right)_\sharp\}.
\]
First, we point out that $\mathcal H_f$ is well defined since any map $g\in W^{1,p}(M',N)$ defined on a compact $m$-dimensional manifold $M'$ induces a homomorphism $g_\sharp:\pi_1(M',\ast)\to\pi_1(N,\ast)$ as follows. Given a generator $\gamma$ for $\pi_1(M',\ast)$, we consider a tubular neighborhood $T\subset M'$ of $\gamma$ in $M$ such that $\psi:\mathbb S^1\times I^{m-1}\to T$ is a smooth immersion, where $I^{m-1}$ is the unit $(m-1)$-cell, and define $\gamma^s:\mathbb S^1\to N$ as $\gamma^s(\cdot):=\psi(\cdot,s)$. Since $M'$ is compact and $p\geq 2$, by H\"older inequality $g\in W^{1,2}(M',N)$ and Proposition 2.3 in \cite{Burstall-London} ensures that there exists $I^{m-1}_g\subseteq I^{m-1}$ such that $I^{m-1}\setminus I^{m-1}_g$ has measure zero and, for all $s,s'\in I^{m-1}_g$, $g$ is continuous on $\gamma^s$ and $g(\gamma^s)$ is homotopic to $g(\gamma^{s'})$. Consequently, for each $\gamma\in \pi_1(M',\ast)$ and $s_0\in I^{m-1}_g$ we can set
\[
g_\sharp[\gamma^{s_0}] = [g(\gamma^{s_0})]
\]
on the generators, and extend $g_\sharp$ so that it is a group homomorphism. By the above considerations, $g_\sharp$ does not depends on the choice of $s_0\in I^{m-1}_g$, while by Proposition 2.4 in \cite{Burstall-London} $g_\sharp$ is also independent of the choice of the generator. \\ 
Since $f\in\mathcal{H}_f$, $\mathcal H_f$ is non-empty and 
\[
\mathcal I_f:=\inf_{v\in\mathcal H_f}E_p(v)<+\infty.
\]
Here $E_p$ denotes the $p$-energy functional defined as
\[
E_p(v):=\frac 1p\int_M |dv|^p dV_M.
\]
Let $\{v_j\}_{j=1}^{\infty}\subset\mathcal H_f$ be a sequence minimizing $p$-energy in $\mathcal H_f$, i.e. $E_p(v_j)\to \mathcal I_f$ as $j\to\infty$. Choosing a subsequence if necessary, we can suppose $E_p(v_j)<2\mathcal I_f$ for all $j$. Fix $k\in\mathbb N$ and consider the sequence $\{v_j|_{M_k}\}_{j=1}^{\infty}$. Let $i:N\hookrightarrow\mathbb R^q$ be an isometric immersion of $N$ into some Euclidean space. Since $i(N)\subset\mathbb R^q$ is compact and $\{E_p(v_j)\}_{j=1}^\infty$ is bounded, $\{v_j|_{M_k}\}_{j=1}^{\infty}$ is bounded in $W^{1,p}(M_k,\mathbb R^q)$ and, up to choosing a subsequence, $v_j|_{M_k}$ converges to some $v^{(k)}\in W^{1,p}(M_k,N)$ weakly in $W^{1,p}$, strongly in $L^p$ and pointwise almost everywhere. This implies $v^{(k)}\in W^{1,p}(M_k,N)$. By the lower semicontinuity of $E_p$ we have
\begin{equation}\label{lsc}
E_p(v^{(k)})\leq \liminf_{j\to\infty}E_p(v_j|_{M_k}).
\end{equation}
We want to show that the homomorphism on fundamental groups induced by $v_j|_{M_k}$ is preserved in the limit $v_j|_{M_k}\to v^{(k)}$. As above, choose a generator $\gamma$ for some fixed class in $\pi_1(M_k,\ast)$ and consider the relative tubular neighborhood $\psi:\mathbb S^1\times I^{m-1}\to T$ and set $I^{m-1}_{v_j}$. For a.e. $s\in I^{m-1}_{v_j}$ there exists a number $K_s$ such that
\begin{equation}\label{eqlim}
\int_{\mathbb S^1}|dv_j(t,s)|^2dt\leq K_s
\end{equation}
for infinitely many $j$. In fact, if by contradiction we assume there exists a set $I'\subset I^{m-1}_{v_j}$ of positive measure such that 
\[
\int_{\mathbb S^1}|dv_j(t,s')|^2dt\to\infty,\qquad\forall s'\in I',
\]
then by H\"older inequality, Fubini's theorem and Fatou's lemma we would have
\begin{align*}
\left(2\mathcal I_f\right)^{\frac2p}\left(\vol \left(\mathbb S^1\times I^{m-1}\right)\right)^{\frac{p-2}p}&\geq
\liminf_{j\to\infty}\int_{\mathbb S^1\times I^{m-1}}|dv_j(t,s)|^2 dV_{\mathbb S^1\times I^{m-1}}\\
&\geq
 \liminf_{j\to\infty}\int_{I'}\int_{\mathbb S^1}|dv_j(t,s)|^2dtds\\
&\geq \int_{I'}\left(\liminf_{j\to\infty}\int_{\mathbb S^1}|dv_j(t,s)|^2dt\right)ds = +\infty.
\end{align*}
So \eqref{eqlim} is proven. The one dimensional Sobolev and Kondrachov's embedding theorems (e.g. \cite{A}, p.53) states that there is a compact immersion $W^{1,2}(\mathbb S^1,N)\hookrightarrow C^0(\mathbb S^1,N)$. By \eqref{eqlim} and the compactness of $N$, $v_j(\cdot,s)$ is uniformly bounded in $W^{1,2}(\mathbb S^1,N)$ so that for a.e. $s\in I^{m-1}$ we get a subsequence $v_j$ converging uniformly on $\gamma^s$. Thus, for a.e. $s\in I^{m-1}$ there is $j$ such that $v_j(\gamma^s)$ is uniformly close to, and hence homotopic to, $v^{(k)}(\gamma^s)$. Hence
\begin{equation}\label{homo_compact}
\left(v^{(k)}\right)_\sharp = \left(f|_{M_k}\right)_\sharp.
\end{equation}
Using standard diagonal arguments we can choose a subsequence of $v_j$ which, for all $k$, converges to $v^{(k)}\in W^{1,p}(M_k,N)$ weakly in $W^{1,p}$, strongly in $L^p$ and pointwise almost everywhere. The map $v_0:M\to N$ which, on $M_k$, takes values $v_0|_{M_k}=v^{(k)}$ is well defined. Indeed, by pointwise convergence, $v^{(k)}$ and $v^{(k+1)}$ agree almost everywhere on $M_k$. Now, $v_0\in W_{loc}^{1,p}(M,N)$ and by \eqref{homo_compact} we get $v_0\in\mathcal{H}_f$. It follows from (\ref{lsc}) and the uniform boundedness of $E_p(v_j|_{M_k})$ that
\begin{align*}
\mathcal I_{f} &\leq  E_p(v_0) = \lim_{k\to\infty} E_p(v_0|_{M_k}) = \lim_{k\to\infty} E_p(v^{(k)}) \\
&\leq \lim_{k\to\infty} \liminf_{j\to\infty}E_p(v_j|_{M_k})\leq  \liminf_{j\to\infty}E_p(v_j) = \mathcal I_f,
\end{align*}
so that $E_p(v_0)=\mathcal I_f$, i.e. $v_0$ minimize the energy in $\mathcal H_{f}$.\\
We are going to show that $v_0\in C^{1,\alpha}$ and, hence, a $p$-harmonic map. To this end, we recall some definitions. A map $v\in W^{1,p}(M,N)$ is said to be $p$\textsl{-minimizing on }$\epsilon$\textsl{-balls} if $E_p(v)\leq E_p(w)$ for any $w\in W^{1,p}(M,N)$ which agrees with $v$ outside some ball $B_r$ of radius $r<\epsilon$, that is, if $v=w$ on $M\setminus B_r$ and $(v-w)|_{B_r}\in W^{1,2}_{0}(B_r,N)$. Moreover, a map $\psi : \mathbb S^{l}\to N$ is said to be a $p$\textsl{-minimizing tangent map} of $\mathbb S^{l}$ if its homogeneous extension $\bar\psi$ to $\mathbb R^{l+1}$ given by 
\[
\bar\psi(x):=\psi\left(\frac x {|x|}\right), \quad\forall x\neq 0,
\]
minimizes the $p$-energy on every compact subset of $\mathbb R^{l+1}$. As observed in \cite{SU} for $p=2$, when restricting to $\mathbb S^{l}$ the $p$-energy of $\bar\psi$ splits in one component tangential to $\mathbb S^{l}$ and one component normal to $\mathbb S^{l}$ which vanishes by homogeneity. Then $\bar\psi$ is $p$-harmonic in $\mathbb R^{l+1}$ if and only if $\psi$ is $p$-harmonic on $\mathbb S^{l}$. Since $\psi$ is a $p$-harmonic map defined on $\mathbb S^{l}$ with values in $N$, which is compact and non-negatively curved, Theorem 1.7 in \cite{Wei-indiana} implies $\psi$ is constant, thus proving that $N$ admits no non-trivial $p$-minimizing tangent maps of $l$ sphere for every $l\geq 1$. On the other hand, choose $\epsilon$ to be less than half the width of the tubular neighborhhoods about the generating curves of various $\pi_1(M,\ast)$ ($\epsilon$ may vary according to the element of $\pi_1(M)$ considered, but this is not important due to the local nature of the regularity results). Then, if $w\in W^{1,p}(M_k,N)$ agrees with $v_0|_{M_k}$ outside some $\epsilon$-ball, we can extend $w$ to $\bar w\in W_{loc}^{1,p}(M,N)$ by setting $\bar w=v_0$ on $M\setminus M_k$, and it is clear that $\bar w\in \mathcal H_f$. Thus
\[
E_p(v_0|_{M_k})+E_p(v_0|_{M\setminus M_k})= E_p(v_0)\leq E_p(\bar w) = E_p(w) + E_p(v_0|_{M\setminus M_k}),
\] 
giving that $v_0$ is $p$-minimizing on $\epsilon$-balls. At this point we can apply a regularity result by Hardt and Lin (see Theorem 4.5 in \cite{HL}) which gives that $v_0$ is $C^{1,\alpha}$ on $M_k$ for each $k$, so $v_0\in C^{1,\alpha}(M,N)$ and is therefore $p$-harmonic since locally $p$-minimizing.\\
It remains to prove that $v_0$ is homotopic to $f$. Since $N$ has non-positive sectional curvatures, $N$ is $K(\pi,1)$, i.e. each homotopy group $\pi_k(N)$ of $N$ is trivial for $k>1$. A standard result says that, in this case, for every compact manifold $M'$ the conjugacy classes of homomorphisms from $\pi_1(M')$ to $\pi_1(N)$ are in bijective correspondence with the homotopy classes of maps from $M'$ to $N$ (see e.g. \cite{Sp} p.428). Thus, the continuous elements of $\mathcal H_f$, and in particular $v_0$, are all homotopic to $f$ on compacta. Finally, to conclude that $v_0$ and $f$ are homotopic as maps from $M$ to $N$, we use the following result attributed to V.L. Hansen; see Theorem 5.1 in \cite{Burstall-London}. We are grateful to professor Fran Burstall for suggesting us the proof Theorem \ref{th_Hansen}.

\begin{theorem}[Hansen]\label{th_Hansen}
Let $M,N$ be connected $C-W$ complexes with $M$ countable and $N$ a $K(\pi,1)$. Let $f,g:M\to N$ be maps that are homotopic on compacta. Then $f,g$ are homotopic as maps from $M$ to $N$. 
\end{theorem}
\end{proof}

We now come to the proof of Theorem \ref{th_SYp}. This combines the techniques introduced by \cite{SY-LN}, \cite{PRS-MathZ} and \cite{HPV}. We consider two non-constant homotopic $p$-harmonic maps $u,v:M\to N$. Since we have no way to compare the two maps in a standard way we proceed as in \cite{SY-LN}, i.e. we consider the map $(u,v):M\to N\times N$ and the distance function on $N$, and their lifting to suitably chosen covering spaces. Nevertheless, the ideas of Schoen and Yau can not be applied directly. Due to the bad behaviour of the composition of $p$-harmonic maps and convex function (by the way, we note that while for $p=2$ the map $(u,v)$ is harmonic, the same does not hold if $p>2$), we are led to introduce special vector field $X$ as in \cite{PRS-MathZ} and apply Kelvin-Nevanlinna-Royden criterion. As in \cite{HPV}, here the criterion has to be applied in its full strength since $\dive X$ has a non-trivial negative part. Indeed, $\dive X$ will depend on the Hessian of the distance function on $N$, but for $p>2$ the Hessian has to be valued in two different direction, which forces to have a negative part. Moreover, since a direct application of the second variation formula for arclength does not suffice, this is obtained through some further computations, see Theorem \ref{lem_hess_p} below. Obviously, throughout all the proof, we have to keep in account the low regularity guaranteed for the $p$-harmonic representatives.

\begin{proof}[Proof (of Theorem \ref{th_SYp})]
We begin by recalling the following interesting characterization of $p$-parabolicity due to \cite{GT} which will be useful later. It goes under the name of Kelvin-Nevanlinna-Royden criterion and was previously proved by T. Lyons and D. Sullivan for $p=2$, \cite{LS}.
 
\begin{theorem}[Gol'dshtein-Troyanov]\label{weak_KNR}
A complete Riemannian manifold $M$ is $p$-parabolic if and only if every vector field $X$ on $M$ such that
\begin{itemize}
\item[(a)] $\left\vert X\right\vert \in L^{\frac{p}{p-1}}\left(  M\right)  $
\item[(b)] $\operatorname{div}X\in L_{loc}^{1}\left(  M\right)  $ and $\min\left(
\operatorname{div}X,0\right)  =\left(  \operatorname{div}X\right)  _{-}\in
L^{1}\left(  M\right)  $
\end{itemize}
satisfies necessarily $0\geq\int_{M}\operatorname{div}X dV_M$.
\end{theorem}

Now, let $u$ and $v$ be two $C^1$ $p$-harmonic maps from $M$ to $N$ which are freely homotopic, and such that $|du|,|dv|\in L^p(M)$. Let $P_M:\tm\to M$ and $P_N:\tn\to N$ be the universal Riemannian covers respectively of $M$ and $N$. Then $\pi_1(M,\ast)$ and $\pi_1(N,\ast)$
act as groups of isometries on $\tm$ and $\tn$ respectively so that $M = \tm/\pi_1(M, \ast)$ and
$N = \tn/\pi_1(N, \ast)$. Let $\dist_{\tn}: \tn \times\tn\to\mathbb R$ be te distance function on $\tn$. Since ${}^{\tn}\sect\leq0$, we know that $\dist_{\tn}$ is smooth on $\tn\times\tn\setminus\tilde D$, where $\tilde D$ is the diagonal set $\{(\tilde x,\tilde x):\tilde x\in\tn\}$, and $\dist_{\tn}^2$ is smooth on $\tn\times\tn$. Now $\pi_1(N,\ast)$ acts on $\tn\times\tn$ as a group of isometries by
\[
\beta(\tilde x,\tilde y) = (\beta(\tilde x),\beta(\tilde y))\qquad\textrm{for }\beta\in\pi_1(N,\ast).
\]
Thus $\dist_{\tn}^2$ induces a smooth function 
\[
\tir^2 : \tilde N_{\times/{}}\to\mathbb R,
\]
where we have defined 
\[
\tilde N_{\times/{}}:=\tn\times\tn / \pi_1(N,\ast).
\]
Let $U: M\times [0, 1]\to N$ be a homotopy
of $u$ with $v$ so that $U(q, 0) = u(q)$ and $U(q, 1) = v(q)$ for all $q\in M$. We choose a lifting $\tilde U: \tm\times [0, 1] \to \tn$, and call $\tilde U(\tq, 0) =: \tilde u(\tq)$ and $\tilde U(\tq, 1) =: \tilde v(\tq)$ for all
$\tq\in \tm$. This defines liftings $\tilde u,\tilde v$ of $u,v$ and, since Riemannian coverings are local isometries, $\tu$ and $\tv$ are $p$-harmonic maps and
\[
|d\tu|(\tq)=|du|(P_M(\tq)),\qquad |d\tv|(\tq)=|dv|(P_M(\tq)).
\]
Now, $\pi_1(M,\ast)$ acts as a group of isometries on $\tm$ and we have 
\begin{equation}\label{gamma}
\tu(\gamma(\tq))=u_\sharp(\gamma)\tu(\tq),\qquad \tv(\gamma(\tq))=v_\sharp(\gamma)\tv(\tq),\qquad\forall\tq\in\tm,\gamma\in\pi_1(M,\ast),
\end{equation}
where $u_\sharp,v_\sharp:\pi_1(M,\ast)\to\pi_1(N,\ast)$ are the induced homomorphism and $u_\sharp\equiv v_\sharp$ since $u$ is homotopic to $v$.\\
Thus, the map $\tilde j: \tm\to\tn\times\tn$ defined by $\tilde j(\tilde x) := (\tilde u(\tilde x), \tilde v(\tilde x))$ induces by (\ref{gamma}) a map
\[
j: M \to \tnn.
\]
Furthermore, we can construct a vector valued $1$-form $J\in T^\ast M\otimes j^{-1}T\tnn$ along $j$ by projecting via \eqref{gamma} the vector valued $1$-form $\tilde J\in T^\ast \tm\otimes \tilde j^{-1}T\left(\tn\times\tn\right)$ along $\tilde j$ defined as
\[
\tilde J := (\Kp\tu,\Kp\tv).
\]   
From now on, the symbol $\Kp\tu$ stands for 
\[
\Kp\tu:=|d\tu|^{p-2}d\tu.
\]
Set $\hat h_A:[0,+\infty)\to\mathbb R$ as $\hat h_A(t):= \sqrt{A+t^2}$ for every $A>1$ and define $h_A:=\hat h_A(\tir)\in C^{\infty}(\tnn,\mathbb R)$. Consider the vector field on $M$ given by
\begin{equation}\label{field_X}
X|_{q}:=\left[dh_A|_{j(q)}\circ J|_{q}\right]^{\sharp}.
\end{equation}
Note that
\begin{equation}\label{field_X'}
X|_{q}:=dP_M|_{\tq} \circ \left.\tilde X\right|_{\tq},
\end{equation}
where
\[
\left.\tilde X\right|_{\tq}:=\left[\left.d\tilde h_A\right|_{\tilde j(\tq)}\circ \left.\tilde J\right|_{\tq}\right]^{\sharp},\qquad \tilde h_A:=\hat h_A\circ\left(\dist_{\tn}^2\right):\tn\times\tn\to\mathbb R.
\]
We claim that \eqref{field_X'} is well defined. To this end, let $S_{\tq}\in T_{\tq}\tm$ be an arbitrary vector and let $\tq'\in P_M^{-1}(q)\subset T\tm$. If $\tq'\neq\tq$, there exists $\gamma\in\pi_1(M,\ast)$ such that $q'=\gamma q$. Then,
\[
\tilde J|_{\gamma\tq}(d\gamma\left( S_{\tq})\right)=\left(d\left[u_{\sharp}(\gamma)\right]\left(\Kp{\tu}(S_{\tq})\right),d\left[v_{\sharp}(\gamma)\right]\left(\Kp{\tv}(S_{\tq})\right)\right).
\]
Since $u$ is homotopic to $v$, $u_{\sharp}=v_{\sharp}$. Moreover $\dist_{\tn}$ is equivariant with respect to the action of $\pi_1(N)$ on $\tn\times\tn$, i.e.
\[
\dist_{\tn}(\beta \tilde x_1, \beta \tilde x_2)=\dist_{\tn}( \tilde x_1, \tilde x_2),\qquad\forall\beta\in\pi_1(N),\ x_1,x_2\in\tn.
\] 
Then, 
\[
dP_M|_{\tq}\circ \left[\left.d\left(\dist_{\tn}^2\right)\right|_{\tilde j(\tq)}\circ \left.\tilde J\right|_{\tq}\right]^{\sharp}
\]
is well defined, and consequently the same holds for $dP_M|_{\tq}\circ\tilde X|_{\tq}$.\\
Now, we want to compute (in the weak sense) $\dive X$ on $M$. We start with the following result, obtained with minor changes from a lemma of Kawai, \cite{Kawai}.
 
\begin{lemma}\label{lem_kawai}
Consider $C^{1}$ $p$-harmonic maps $u,v:M\to N$ and a smooth function $h:N\times N\to\mathbb R$. Then the identity
\begin{align}\label{weak_kawai}
{}^M\tr &{}^{N\times N}\Hess h|_{(u,v)}\left(\left(du,dv\right),\left(\Kp u,\Kp v\right)\right) \\
&= {}^M\dive \left[dh|_{(u,v)}\circ \left(\Kp u,\Kp v\right)\right]^\sharp,\nonumber
\end{align}
holds weakly on $M$.
\end{lemma}
\begin{proof}
Consider $(u,v):M\to N\times N$. Let $\eta\in C^{\infty}(M,\mathbb R)$ be a compactly supported function and define a vector valued $1$-form $\psi\in T^\ast M\otimes (u,v)^{-1}T (N\times N)$ along $(u,v)$ as $\psi:=D(\eta\nabla h|_{(u,v)})$, that is
\[
\psi(V) = (d\eta(V))\;{}^{N\times N}\nabla h|_{(u,v)} + \eta\; {}^{N\times N}\nabla_{d(u,v)(V)}{}^{N\times N}\nabla h|_{(u,v)}
\]
for all vector fields $V$ on $M$. Since $u,v$ are $p$-harmonic, by the structure of Riemannian products we have that $\dive\left(\Kp u,\Kp v\right)=0$ weakly on $M$, that is
\[
\int_M\left\langle \xi, \left(\Kp u,\Kp v\right)\right\rangle_{HS_{\times}}=0,\quad\forall \xi\in T^\ast M\otimes (u,v)^{-1}T \left(N\times N\right),
\]
where $\left\langle ,\right\rangle_{HS_{\times}}$ is the Hilbert-Schmidt scalar product on $T^\ast M\otimes (u,v)^{-1}T \left(N\times N\right)$. Hence, choosing $\xi=\psi$ in this latter, we obtain
\begin{align*}
0&=\int_M\left\langle \psi, \left(\Kp u,\Kp v\right)\right\rangle_{HS_{\times}}\\
&=\int_M\left\langle d\eta(\cdot)\otimes{}^{N\times N}\nabla h|_{(u,v)}, \left(\Kp u,\Kp v\right)\right\rangle_{HS_{\times}}\\
&+\int_M\left\langle \eta\; {}^{N\times N}\nabla_{d(u,v)(\cdot)}{}^{N\times N}\nabla h|_{(u,v)}, \left(\Kp u,\Kp v\right)\right\rangle_{HS_{\times}}\\
&=\int_M \left[dh|_{(u,v)}\circ \left(\Kp u,\Kp v\right)\right]({}^M\nabla\eta)\\
&+\int_M \eta\;{}^M\tr\left\langle {}^{N\times N}\nabla_{(du,dv)}{}^{N\times N}\nabla h|_{(u,v)},\left(\Kp u,\Kp v\right)\right\rangle_{N\times N},
\end{align*}
which turns to be the weak formulation of \eqref{weak_kawai}.
\end{proof}

According to Lemma \ref{lem_kawai}, because of \eqref{field_X} and since $\pi_1(M,\ast)$ acts on $\tm$ as a group of isometries, we have that for all $q\in M$ and for any choice of $\tq\in P_M^{-1}(q)$
\begin{align}\label{div_hes}
{}^M\dive X|_q &= {}^{\tm}\dive\tilde X|_{\tq} = {}^{\tm}\dive\left[d\tilde h_A|_{\tilde j(\tq)}\circ \tilde J\right]^{\sharp}\\
&={}^{\tm}\tr{}^{\tn\times \tn}\Hess \tilde h_A|_{\tilde j(\tq)}\left(d\tilde j,\tilde J\right) = {}^{M}\tr {}^{\tnn}\Hess h_A|_{j(q)}\left(dj,J\right)\nonumber
\end{align}
holds weakly on $M$. Observe that
\begin{equation}\label{dh_A}
dh_A=\frac{d\frac{\tir^2}2}{\sqrt{A+\tir^2}} = \frac{\tir d\tir}{\sqrt{A+\tir^2}} 
\end{equation}
and
\[
{}^{\tnn}\Hess h_A = \frac{{}^{\tnn}\Hess \tir^2}{2\sqrt{A+\tir^2}} - \frac{\tir^2}{(A+\tir^2)^{3/2}}d\tir\otimes d\tir.
\]
Then, in order to deal with $\dive X$, we want to compute 
\begin{equation}\label{2hess}
{}^{\tnn}\Hess \tir^2|_{j(q)}\left(dj,J\right)={}^{\tn\times \tn}\Hess \dist_{\tn}^2|_{\tilde j(\tq)}\left(d\tilde j,\tilde J\right),\qquad\tq\in P_M^{-1}(q).
\end{equation}

\begin{theorem}\label{lem_hess_p}
Suppose $N$ is a simply connected Riemannian manifolds such that $\nsect\leq 0$ and fix points $u,v$ in $N$. Let ${}^N r: N\times N\to\mathbb R$ be defined by $\nr u v :={}^N\operatorname{dist}(u,v)$ and let $X=X_1+X_2\in T_{(u,v)}N\times N$, with $X_1\in T_uN$ and $X_2\in T_vN$. Then , for every $p\geq 2$,
\[
\left.^{N\times N}\Hess{}^Nr^2\right|_{(u,v)}(X,(|X_1|^{p-2}X_1,|X_2|^{p-2}X_2))\geq 0
\]
and the equality holds if and only if there is a parallel vector field $Z$, defined along the unique geodesic $\gamma_{uv}$ joining $u$ and $v$, such that $Z(u)=X_1$, $Z(v)=X_2$ and $\nmetr{{}^N R(Z,T)T}Z\equiv0$ along $\gamma_{uv}$. Moreover, $d({}^Nr)(X)=0$.\\
In particular, if $\nsect<0$, $Z$ is proportional to $T$.
\end{theorem}

We begin with the following lemma

\begin{lemma}\label{milnor}
Consider a Riemannian manifold $Q$ and a function $f\in C^2(Q,\mathbb R)$. Let $q\in Q$ and $X,Y\in T_qQ$. Let $\eta_X:[-\epsilon,\epsilon]\to Q$ be the constant speed geodesic s.t. $\eta_X(0)=q$ and $\dot\eta_X(0)=X$. Moreover, define $Y_s\in T_{\eta_X(s)}Q$ as the vectors obtained by parallel translating $Y=Y_0$ along $\eta_X$, and let $\eta_Y^{(s)}:[-\delta,\delta]\to Q$ be the constant speed geodesic s.t. $\eta_Y^{(s)}(0)=\eta_X(s)$ and $\dot\eta_Y^{(s)}(0)=Y_s$. Then
\[
{}^Q\Hess f|_q(X,Y)=\left.\frac{\partial^2}{\partial s\partial t}\right|_{s=t=0}f(\eta_Y^{(s)}(t)).
\]  
\end{lemma}

\begin{proof}
We have
\begin{align*}
\left.\frac{\partial^2}{\partial s\partial t}\right|_{s=t=0}f(\eta_Y^{(s)}(t)) &= \left.\frac{\partial}{\partial s}\right|_{s=0}\left.\left\langle \nabla f, \dot\eta_Y^{(s)}\right\rangle\right|_{t=0}
=\left.\frac{\partial}{\partial s}\right|_{s=0}\left\langle \nabla f|_{\eta_X(s)}, Y_s\right\rangle\\ 
&=\left.\left\langle \nabla_{\dot\eta_X(s)}\nabla f, Y_s\right\rangle\right|_{s=0} + \left.\left\langle \nabla f, \nabla_{\dot\eta_X(s)}Y_s\right\rangle\right|_{s=0}\\
&=\left\langle \nabla_{X}\nabla f, Y\right\rangle={}^Q\Hess f|_q(X,Y),
\end{align*}
since $\nabla_{\dot\eta_X(s)}Y_s\equiv 0$ by construction.
\end{proof}

\begin{proof}[Proof (of Theorem \ref{lem_hess_p})]
As above, for the ease of notation, for each vector field $\xi$ we set $k_p(\xi):=|\xi|^{p-2}\xi$. Define the vector field $Y\in T_{(u,v)}N\times N$ as $Y=(Y_1,Y_2)=(k_p(X_1),k_p(X_2))$. Let $D$ be the diagonal set 
\[
D:=\{ (u_1,u_1):u_1\in N \}\subset N\times N 
\]
so that ${}^N r$ is smooth on $N\times N\setminus D$, ${}^N r^2$ is smooth on $N\times N$ and for every pair $(u,v)\in N\times N \setminus D$ there is a unique shortest geodesic from $u$ to $v$. We call $\gamma_{u,v}$ such a geodesic parametrized by arc length.\\
Let $\sigma^X:[-\epsilon,\epsilon]\to N\times N$, be the constant speed geodesic on $N\times N$ satisfying $\sigma^X(0)=(u,v)$ and $\dot\sigma^X(0)=X=(X_1,X_2)$. We then have $\sigma^X=(\sigma^X_1,\sigma^X_2)$ where $\sigma^X_1$ and $\sigma^X_2$ are geodesic on $N$ satisfying $\sigma^X_1(0)=u$, $\sigma^X_2(0)=v$ and $\dot{\sigma^X_i}(0)=X_i$, $i=1,2$. As in Lemma \ref{milnor}, let $Y_s$ be the vector field along $\sigma^X$ obtained by parallel transport of $Y$ and let $\sigma^{(s),Y}:[-\delta,\delta]\to N\times N$ be the constant speed geodesic s.t. $\sigma^{(s),Y}(0)=\sigma^X(s)$ and $\dot\sigma^{(s),Y}(0)=Y_s$. As above we can split $\sigma^{(s),Y}$ in two geodesic of $N$, i.e. $\sigma^{(s),Y}=(\sigma^{(s),Y}_1,\sigma^{(s),Y}_2)$. \\
Set $\bar R:={}^N r(u,v)$ and, for each couple of points $y_1,y_2\in N$ let $\gamma_{y_1,y_2}:[0,\bar R]\to N$ be the (unique) constant speed geodesic joining $y_1$ and $y_2$.\\
At this point, we can consider a two parameters geodesic variation of $\gamma_{u,v}$ defining $\alpha:[0,\bar R]\times [-\epsilon,\epsilon]\times [-\delta,\delta]\to N$ as
\[
\alpha (t,z,w) := \gamma_{\sigma^{(z),Y}_1(w),\sigma^{(z),Y}_2(w)}(t).
\]
We now define the variational vector fields
\begin{align*}
&\hat Z(t,z,w):=\frac{\partial}{\partial z}\alpha(t,z,w),\quad \hat W(t,z,w):=\frac{\partial}{\partial w}\alpha(t,z,w),\\
&Z(t):=\hat Z(t,0,0),\quad W(t):=\hat W(t,0,0),\quad T(t):=\frac{\partial}{\partial t}\alpha(t,0,0)=\dot\gamma_{u,v}(t).
\end{align*}
Here, we are using the notation
\[
\frac{\partial}{\partial z}\alpha(t,z,w) := \left.d\alpha\right|_{(t,z,w)} \left(\frac{\partial}{\partial z}\right), 
\]
with $\frac{\partial}{\partial z}=\left(0,\frac{\partial}{\partial z},0\right)\in T\left([0,\bar R]\times [-\epsilon,\epsilon]\times [-\delta,\delta]\right)$.
Since both the one parameter variations $\alpha(t,z,0)$ and $\alpha(t,0,w)$ are geodesic variations, we have that $Z$ and $W$ are the corresponding Jacobi fields along $\gamma_{u,v}$. Then they satisfy
\[
\begin{array}{ll}
Z(0)=X_1, &\quad W(0)=Y_1= \kp{X_1}= \kp{Z(0)}, \\ Z(\bar R)=X_2, &\quad W(\bar R)=Y_2=\kp{X_2}= \kp{Z(\bar R)},
\end{array}
\]
and the Jacobi equations
\begin{align*}
\nabla_T\nabla_T Z+ {}^NR(Z,T)T = 0 = \nabla_T\nabla_T W+ {}^NR(W,T)T.
\end{align*}
For each $z\in[-\epsilon,\epsilon]$ and $w\in[-\delta,\delta]$, let 
\[
L_{\alpha}(z,w):=\int_0^{\bar R}\left|\frac{\partial}{\partial t}\alpha(t,z,w)\right| dt 
\]
be the length of the geodesic curve $t\mapsto \alpha (t,z,w)$. By Lemma \ref{milnor} we have
\begin{align*}
{}^{N\times N}\Hess{}^Nr|_{(u,v)}(X,Y) &= \left.\frac{\partial^2}{\partial z\partial w}\right|_{z=w=0}{}^Nr(\sigma^{(z),Y}(w))\\
&=\left.\frac{\partial^2}{\partial z\partial w}\right|_{z=w=0}{}^Nr(\sigma_1^{(z),Y}(w),\sigma_2^{(z),Y}(w))\\
&=\left.\frac{\partial^2}{\partial z\partial w}\right|_{z=w=0} L_{\alpha}(z,w).
\end{align*}
On the other hand, by the second variation of arc length (see \cite{CE}, page 20) we have
\begin{align}\label{2nd_var}
\left.\frac{\partial^2}{\partial z\partial w}\right|_{z=w=0} L_{\alpha}(z,w) &= \left.\nmetr{\nabla_{\hat Z}\hat W(t,0,0)}{T(t)}\right|_{t=0}^{t=\bar R}+ \int_0^{\bar R} \nmetr{\nabla_TZ}{\nabla_TW}\\ 
&-\int_0^{\bar R}\nmetr{{}^NR(W,T)T}Z - \int_0^{\bar R}T\nmetr ZTT\nmetr WT .\nonumber
\end{align}
We note that the vector fields $\hat Z$ and $\hat W$ are defined along the map $\alpha$. Accordingly, the covariant derivative at the first term on RHS of \eqref{2nd_var} has the meaning, 
\[
\nabla_{\hat Z}\hat W(t,0,0)=\left.\nabla_{\hat Z}\right|_{\alpha(t,0,0)}\left(\frac{\partial}{\partial w}\alpha(t,z,0)\right).
\]
First, observe that, by construction of $\alpha$ and due to the choice of the geodesics $\sigma^{(z),Y}_1(w)$ and $\sigma^{(z),Y}_2(w)$, we have
\[
\nabla_{\hat Z}\hat W(0,0,0)=\nabla_{\hat Z}\hat W(\bar R,0,0)=0,
\]
which implies that the first term on RHS of \eqref{2nd_var} vanishes. Moreover, using the Jacobi equation for $Z$ and the values of $W$ at $t=0$ and $t=\bar R$ we can compute
\begin{align}\label{1st}
&\int_0^{\bar R} \left\{\nmetr{\nabla_TZ}{\nabla_TW}-\nmetr{{}^NR(W,T)T}Z\right\}\\
&=\int_0^{\bar R} \left\{\nmetr{\nabla_TZ}{\nabla_TW}+\nmetr{\nabla_T\nabla_TZ}W\right\}\nonumber\\
&=\int_0^{\bar R} T\nmetr{\nabla_TZ}W \nonumber\\
&= \left.\nmetr{\nabla_TZ}W\right|_{t=0}^{t=\bar R}\nonumber\\
&=\left.\nmetr{\nabla_TZ}{\kp Z}\right|_{t=0}^{t=\bar R}\nonumber\\
&=\int_0^{\bar R} T\nmetr{\nabla_TZ}{\kp Z} \nonumber\\
&=\int_0^{\bar R} \left\{\nmetr{\nabla_T\nabla_T Z}{\kp Z}+ T\left(|Z|^{p-2}\right)\nmetr{\nabla_T Z}Z+|Z|^{p-2}|\nabla_TZ|^2\right\}\nonumber\\
&=\int_0^{\bar R} \left\{-|Z|^{p-2}\nmetr{^{N}R(Z,T)T}{Z}+\frac{1}{2}T\left(|Z|^2\right)T\left(|Z|^{p-2}\right)+|Z|^{p-2}|\nabla_TZ|^2\right\}. \nonumber
\end{align}
Since $T$ is parallel, the Jacobi equation implies
\begin{align}\label{T_par}
TT\nmetr ZT &= T\nmetr{\nabla_TZ}T \\
&= \nmetr{\nabla_T\nabla_T Z}T \nonumber\\
&= \nmetr{{}^NR(T,Z)T}T = 0.\nonumber
\end{align}
Then
\begin{align}\label{2nd}
&\int_0^{\bar R}T\nmetr ZTT\nmetr WT \\
&= \int_0^{\bar R}\left\{T\left(\nmetr WTT\nmetr ZT\right)- \nmetr WTTT\nmetr ZT\right\} \nonumber\\
&=\left.\left(\nmetr WTT\nmetr ZT\right)\right|_{t=0}^{t=\bar R}\nonumber\\
&=\left.\left(\nmetr {\kp Z}TT\nmetr ZT\right)\right|_{t=0}^{t=\bar R}\nonumber\\
&=\int_0^{\bar R}T\left(\nmetr {\kp Z}TT\nmetr ZT\right)\nonumber\\
&=\int_0^{\bar R} \left\{T\left(|Z|^{p-2}\right)\nmetr ZT T\nmetr ZT + |Z|^{p-2}\left(T\nmetr ZT\right)^2\right\}.\nonumber
\end{align}
Inserting \eqref{1st} and \eqref{2nd} in \eqref{2nd_var} we get
\begin{align}\label{res}
&{}^{N\times N}\Hess{}^Nr|_{(u,v)}(X,Y)\\
&=\int_0^{\bar R}|Z|^{p-2}\left\{-\nmetr{{}^NR(Z,T)T}Z+|\nabla_TZ|^2-\left(T\nmetr ZT\right)^2 \right\}\nonumber\\
&+\int_0^{\bar R}\frac12 T\left(|Z|^{p-2}\right)T\left(|Z|^{2}\right)-\int_0^{\bar R} T\left(|Z|^{p-2}\right)\nmetr ZT T\nmetr ZT.\nonumber
\end{align}
We consider the three integrals separately. First, since $\left|T\nmetr ZT\right|=\left|\nabla_TZ^{T}\right|$, we have that the first integral at RHS of \eqref{res} is equal to
\begin{align}\label{1res}
\int_0^{\bar R}|Z|^{p-2}\left\{|\nabla_TZ^{\perp}|^2-\nmetr{{}^NR(Z,T)T}Z\right\},
\end{align}
where $Z^T$ and $Z^{\perp}$ denote the components of $Z$ respectively parallel and normal to $T$, and the integral is positive by the curvature assumptions on $N$.\\
As for the second integral, assume $2<p<4$, the other cases being easier. We have
\begin{align}\label{2res}
T\left(|Z|^{p-2}\right)T\left(|Z|^{2}\right)=\frac2{p-2}|Z|^{4-p}\left[T\left(|Z|^{p-2}\right)\right]^2\geq0.
\end{align}
Finally, recall \eqref{T_par} and note that this implies that $T\nmetr ZT$ is constant along $\gamma_{u,v}$ and takes value 
\begin{align}\label{dr1}
T\nmetr ZT \equiv \frac 1 {\bar R}\left(\left.\nmetr ZT\right|_{t=\bar R} - \left.\nmetr ZT\right|_{t=0}\right).
\end{align}
On the other hand, we have that
\begin{align}\label{dr2}
\left.d{}^Nr\right|_{(u,v)}(X)
&=\left.d{}^Nr\right|_{(u,v)}((X_1,X_2))\\
&=\left.d{}r_v\right|_{u}(X_1) + \left.dr_u\right|_{v}(X_2)\nonumber\\
&=-\nmetr{X_1}{\dot\gamma_{u,v}(0)}+\nmetr{X_2}{\dot\gamma_{u,v}(\bar R)},\nonumber
\end{align}
where $r_u,r_v:N\to\mathbb R$ are defined as $r_u(\cdot):={}^Nr(u,\cdot)$ and $r_v(\cdot):={}^Nr(\cdot,v)$. Combining \eqref{dr1} and \eqref{dr2} we get
\begin{align}\label{dr3}
T\nmetr ZT \equiv \frac{d{}^Nr(X)}{\bar R},
\end{align}
which in turn implies
\begin{align}\label{3res}
&\int_0^{\bar R} T\left(|Z|^{p-2}\right)\nmetr ZT T\nmetr ZT\\ 
&=\int_0^{\bar R} T\left(|Z|^{p-2}\nmetr ZT T\nmetr ZT\right) - \int_0^{\bar R} |Z|^{p-2}\left(T\nmetr ZT\right)^2\nonumber\\
&=\frac{d{}^Nr(X)}{\bar R}\left[|Z|^{p-2}\nmetr ZT\right]_{t=0}^{t=\bar R} - \left(\frac{d{}^Nr(X)}{\bar R}\right)^2\int_0^{\bar R}|Z|^{p-2}.\nonumber
\end{align}
Moreover, reasoning as for \eqref{dr3}, we compute
\begin{align}\label{4res}
&\frac{d{}^Nr(X)}{\bar R}\left[|Z|^{p-2}\nmetr ZT\right]_{t=0}^{t=\bar R} \\
&= \frac{d{}^Nr(X)}{\bar R}\left[\nmetr {\kp{X_2}}{\dot\gamma_{u,v}(\bar R)}-\nmetr {\kp{X_1}}{\dot\gamma_{u,v}(0)}\right] \nonumber\\
&= \frac{d{}^Nr(X)d{}^Nr(Y)}{\bar R}\nonumber.
\end{align}
Combining \eqref{res}, \eqref{1res}, \eqref{3res} and \eqref{4res}, we obtain
\begin{align*}
&{}^{N\times N}\Hess{}^Nr|_{(u,v)}(X,Y)\\
&=
\int_0^{\bar R}|Z|^{p-2}\left\{|\nabla_TZ^{\perp}|^2-\nmetr{{}^NR(Z,T)T}Z\right\}+\frac12\int_0^{\bar R} T\left(|Z|^{p-2}\right)T\left(|Z|^{2}\right)\nonumber\\
&-\frac{d{}^Nr(X)d{}^Nr(Y)}{\bar R} + \left(\frac{d{}^Nr(X)}{\bar R}\right)^2\int_0^{\bar R}|Z|^{p-2}.\nonumber
\end{align*}
Finally, since 
\[
\Hess r^2 = 2r\Hess r + 2dr\otimes dr,
\]
recalling also \eqref{2res}, we get
\begin{align}\label{res2}
&{}^{N\times N}\Hess{}^Nr^2|_{(u,v)}(X,Y)\\
&=
2\bar R\int_0^{\bar R}|Z|^{p-2}\left\{|\nabla_TZ^{\perp}|^2-\nmetr{{}^NR(Z,T)T}Z\right\}\nonumber\\
&+\bar R\int_0^{\bar R} T\left(|Z|^{p-2}\right)T\left(|Z|^{2}\right)+2\frac{\left(d{}^Nr(X)\right)^2}{\bar R}\int_0^{\bar R}|Z|^{p-2}\geq 0.\nonumber
\end{align}
This conclude the first part of the proof. Now, assume 
\[
{}^{N\times N}\Hess{}^Nr^2|_{(u,v)}(X,Y)=0.
\]
From \eqref{res2} we get that $d{}^Nr(X)=0$, $\nmetr{{}^NR(Z,T)T}Z\equiv 0$ along $\gamma_{u,v}$ and, using also \eqref{T_par},
\[
|\nabla_TZ|^2=|\nabla_TZ^{\perp}|^2+|\nabla_TZ^T|^2\equiv 0,
\]
that is $Z$ is parallel along $\gamma_{u,v}$. 
\end{proof}

Recalling \eqref{2hess} and applying Theorem \ref{lem_hess_p} with $N=\tn$ and, with an abuse of notation, $X=d\tilde j$ we get

\begin{corollary}\label{coro_hess_p}
With the definitions introduced above, for all $q\in M$, $E_q\in T_qM$ and for any choice of $\tq\in P_M^{-1}(q)$ and $\tilde E_q = \left[d(P_M)|_{\tq}\right]^{-1}(E_q)$ we have
\[
\left.^{\tnn}\Hess\tir^2\right|_{j(q)}(dj|_q(E_q),\left.J\right|_q(E_q))\geq 0
\]
and the equality holds if and only if there is a parallel vector field $Z$, defined along the unique geodesic $\tilde\gamma_{\tq}$ in $\tn$ joining $\tilde u(\tq)$ and $\tilde v(\tq)$, such that $Z(\tu(\tq))=d\tu|_{\tq}(\tilde E_q)$, $Z(\tv(\tq))=d\tv|_{\tq}(\tilde E_q)$ and $\nmetr{{}^N R(Z,\dot{\tilde\gamma}_{\tq})\dot{\tilde\gamma}_{\tq}}Z\equiv0$ along $\tilde\gamma_{\tq}$. Moreover, $d(\operatorname{\dist}_{\tn})(d\tilde j(\tilde E_q))=0$.\\
In particular, if $\nsect<0$, $Z$ is proportional to $\dot{\tilde\gamma}_{\tq}$.
\end{corollary}

We go back to the proof of Theorem \ref{th_SYp}. From \eqref{div_hes}, \eqref{2hess}, applying Corollary \ref{coro_hess_p} and observing that
\[
\frac{t^2}{(A+t^2)^{3/2}}\leq \frac 1 {\sqrt{A+t^2}}\leq A^{-1/2},\quad\forall t>0,
\]
we get
\begin{align}\label{res3}
{}^M\dive X|_q &= \frac {{}^{M}\tr {}^{\tnn}\Hess \tir^2|_{j(q)}\left(dj_{q},J|_{q}\right)}{2\sqrt{A+\tir^2(j(q))}}\\
&- \frac{\tir^2(j(q))}{(A+\tir^2(j(q)))^{3/2}} {}^M\tr \left[d\tir|_{j(q)}\left(dj|_{q}\right)  \left.d\tir\right|_{j(q)}\left(J|_{q}\right)\right]\nonumber\\
&\geq -A^{-1/2}\left(|du|(q)+|dv|(q)\right)\left(|du|^{p-1}(q)+|dv|^{p-1}(q)\right)\nonumber\\
&\geq -2A^{-1/2}\left(|du|^{p}(q)+|dv|^{p}(q)\right)\nonumber,
\end{align}
from which
\begin{align}\label{res4}
\left(\dive X|_q\right)_-\leq 2A^{-1/2}\left(|du|^{p}(q)+|dv|^{p}(q)\right)\in L^1(M).
\end{align}
Moreover, since $t/\sqrt{A+t^2}<1$, \eqref{dh_A} implies
\begin{align}\label{res5}
|X|^{\frac{p}{p-1}}(q) = &\leq \left(|du|^{p-1}(q)+|dv|^{p-1}(q)\right)^{\frac{p}{p-1}}\\
&\leq 2^{\frac{1}{p-1}}\left(|du|^{p}(q)+|dv|^{p}(q)\right)\in L^1(M).\nonumber
\end{align}
For every $T>0$, set
\[
M_T = \{q\in M : \tir(j(q))\leq T\}\qquad\textrm{and}\qquad M^T:=M\setminus M_T.
\]
From \eqref{res4} and \eqref{res5}, we can apply Proposition \ref{weak_KNR} to deduce that 
\[
\int_M{}^M\dive X\leq 0,
\]
which by \eqref{res3} gives
\begin{align}\label{res6}
&\int_M{\frac{\tir^2(j(q))}{(A+\tir^2(j(q)))^{3/2}}{}^M\tr\left[d\tir|_{j(q)}\left(dj|_{q}\right)  \left.d\tir\right|_{j(q)}\left(J|_{q}\right)\right]}\\
&\geq
\int_M \frac {{}^{M}\tr {}^{\tnn}\Hess \tir^2|_{j(q)}\left(dj|_{q},J|_{q}\right)}{2\sqrt{A+\tir^2(j(q))}}
\nonumber\\ &
\geq \int_{M_T} \frac {{}^{M}\tr {}^{\tnn}\Hess \tir^2|_{j(q)}\left(dj|_{q},J|_{q}\right)}{2\sqrt{A+\tir^2(j(q))}}\nonumber\\ 
&\geq \frac{1}{2\sqrt{A+T}} \int_{M_T} {}^{M}\tr {}^{\tnn}\Hess \tir^2|_{j(q)}\left(dj|_{q},J|_{q}\right)
\geq 0. \nonumber
\end{align}
The real valued function $t\mapsto\frac{t}{(A+t)^{3/2}}$ has a
global maximum at $t=2A$, is increasing in $(0,2A)$ and satisfies 
\[
\frac{t}{(A+t)^{3/2}}<\frac{1}{(A+t)^{1/2}}.
\]
Hence, up
to choosing $A>T^2/2$, we have 
\begin{align}\label{res7}
&\int_M{\frac{\tir^2(j(q))}{(A+\tir^2(j(q)))^{3/2}}{}^M\tr\left[d\tir|_{j(q)}\left(dj|_{q}\right)  \left.d\tir\right|_{j(q)}\left(J|_{q}\right)\right]}\\
&\leq\frac{T^2}{(A+T^2)^{3/2}} \int_{M_T}2\left(|du|^p+|dv|^p\right) + \frac{1}{\sqrt{A+T^2}} \int_{M^T}2\left(|du|^p+|dv|^p\right).\nonumber
\end{align}
Inserting \eqref{res7} in \eqref{res6} we get
\begin{align*}
\int_{M_T} {}^{M}\tr {}^{\tnn}\Hess \tir^2|_{j(q)}\left(dj|_{q},J|_{q}\right)
&\leq \frac{4T^2}{A+T^2} \int_{M_T}\left(|du|^p+|dv|^p\right) \\
&+ 4\int_{M^T}\left(|du|^p+|dv|^p\right)\nonumber,
\end{align*}
and letting $A\to+\infty$ this latter gives
\begin{align*}
\int_{M_T} {}^{M}\tr {}^{\tnn}\Hess \tir^2|_{j(q)}\left(dj|_{q},J|_{q}\right)\leq 4\int_{M^T}\left(|du|^p+|dv|^p\right)\nonumber.
\end{align*}
Since $|du|,|dv|\in L^p(M)$ we can let $T\to+\infty$, applying respectively monotone and dominated convergence to LHS and RHS integrals, thus obtaining
\begin{equation}\label{null}
\int_{M} {}^{M}\tr {}^{\tnn}\Hess \tir^2|_{j(q)}\left(dj|_{q},J|_{q}\right)=0.
\end{equation}
Fix an orthonormal frame $\{E_i\}_{i=1}^m$ for $M$. Then \eqref{null} gives
\[
{}^{\tnn}\Hess \tir^2|_{j(q)}\left(dj(E_i),J(E_i)\right)=0,
\]
for all $i=1,\dots,m$ and $\tq\in M$. At this point, applying again Corollary \ref{coro_hess_p} implies
\[
d\left(\operatorname{dist}_{\tn}\right)\left(d\tu(\tilde E_i),d\tv(\tilde E_i)\right)=d\left(\operatorname{dist}_{\tn}\circ(\tu,\tv)\right)(\tilde E_i)\equiv 0
\]
and, since $\{\tilde E_i\}_{i=1}^m$ span all $T_{\tq}\tm$, we get that $\left(\operatorname{dist}_{\tn}\circ(\tu,\tv)\right)$ is constant on $\tm$. Accordingly, for each $\tq\in\tm$ the unique geodesic $\tilde\gamma_{\tq}$ from $\tu(\tq)$ to $\tv(\tq)$ can be parametrized on $[0,1]$ proportional (independent of $\tq$) to arclength. We define a one-parameter family of maps $\tu_t:\tm\to\tn$ by letting $\tu_t(\tq):=\tilde\gamma_{\tq}(t)$. Then we see that $\tu_0=\tu$ and $\tu_1=\tv$. Corollary \ref{coro_hess_p} states also that for each $i=1,\dots,m$ there exists a parallel vector field $Z_i$, defined along $\tilde\gamma_{\tq}$ in $\tn$, such that $Z_i(0)=d\tu|_{\tq}(\tilde E_i)$, $Z_i(1)=d\tv|_{\tq}(\tilde E_i)$ and $\nmetr{{}^N R(Z_i,\dot{\tilde\gamma}_{\tq})\dot{\tilde\gamma}_{\tq}}{Z_i}\equiv0$ along $\tilde\gamma_{\tq}$. In particular $Z$ is a Jacobi field along $\tilde\gamma_{\tq}$. By the proof of Theorem \ref{lem_hess_p} it turns out that 
\begin{equation}\label{z-tu}
Z_i(t)\equiv \left.d\tu_t\right|_{\tq}(\tilde E_i).
\end{equation}
In fact, let $\zeta_i:(-\varepsilon,\varepsilon)\to \tm$, $\varepsilon>0$, be a smooth curve such that $\dot\zeta(0)=\tilde E_i$. By definition of differential we have that 
\[
d\tu_t(\tilde E_i)=\left.\frac{\partial}{\partial s}\right|_{s=0}\left(\tu_t\circ\zeta\right)(s).
\]
On the other hand, since $\left(\tu_t\circ\zeta\right) (s)=\gamma_{\zeta(s)}(t)$, we get that $d\tu_t(\tilde E_i)$ is the variational field of the geodesic variation
\[
(t,s)\mapsto \gamma_{\zeta(s)}(t),
\]
then $d\tu_t(\tilde E_i)$ is a Jacobi field along $\tilde\gamma_{\tq}$ and, by the uniqueness of the Jacobi fields with given boundary values, \eqref{z-tu} is proved.\\
In the special situation $\nsect<0$, for all $\tq\in\tm$ and $i=1,\dots,m$, the parallel vector field $Z_i$ along $\gamma_{\tq}$ has to be proportional to $\dot\gamma_{\tq}$. Hence $\tu(\tm)$ and $\tv(\tm)$ have to be contained in a geodesic of $\tn$ and projecting on $M$ we get the proof of case $i)$ of Theorem \ref{th_SYp}.\\
In general, because of the equivariance property \eqref{gamma} and by the uniqueness of the construction above, for all $\gamma\in\pi_1(M,\ast)$ and $t\in [0,1]$ we have that 
\begin{align}\label{betagamma}
\tu_t\circ\gamma= \beta\circ\tu_t,
\end{align}
where $\beta=u_{\sharp}(\gamma)=v_{\sharp}(\gamma)\in \pi_1(N,\ast)$. Thus we have induced maps $u_t:M\to N$ for $t\in [0,1]$ such that $u_0\equiv u$ and $u_1\equiv v$.\\
Let $\gamma_q(t)$ be the geodesic from $u(q)$ to $v(q)$ in $M$ obtained by projection from $\tilde\gamma_{\tq}$. Projecting $Z_i$, which is equivariant by \eqref{z-tu} and \eqref{betagamma}, the identity \eqref{z-tu} implies $du_t$ is a parallel vector field along $\gamma_q$. Therefore, the $p$-energy density of $u_t$
\[
e_p(u_t)(q):=\left(\sum_{i=1}^m\left|du_t(E_i)\right|^2\right)^{\frac{p}2}
\]
is constant along $\gamma_q$ for each $q\in M$ and, consequently, the $p$-energies of $u_t$ satisfy
\begin{equation}\label{energies}
E_p(u)=E_p(u_t)=E_p(v),\qquad\forall t\in[0,1],
\end{equation}
that is, every $p$-harmonic map of finite $p$-energy homotopic to $u$ has the same $p$-energy as $u$.\\
Now, suppose $N$ is compact. In case also $M$ is compact, Corollary 7.2 in \cite{Wei-indiana} immediately implies $u_t$ is $p$-harmonic for all $t\in[0,1]$. Otherwise, by Theorem \ref{th_Wei-existence} we know that there exists a $p$-harmonic map $u_{t,\infty}\in \mathcal H_{u_t}$ which minimizes $p$-energy in the homotopy class of $u_t$, which, by construction, is the same homotopy class of $u$. Applying \eqref{energies} with $v=u_{t,\infty}$ we have
\begin{equation*}\label{energies2}
E_p(u_{t,\infty})=E_p(u)=E_p(u_t).
\end{equation*}
On the other hand, if we assume that $u_t$ is not $p$-harmonic, for each $\epsilon>0$ there exists an $\epsilon$-ball $B_{\epsilon}$ such that $u_t$ does not minimize energy on $B_\epsilon$. Namely, there exists a map $\hat u_{t,\epsilon}$ such that 
\[
E_p(\hat u_{t,\epsilon}|_{B_\epsilon})<E_p(u_{t}|_{B_\epsilon}),
\]
so that extending $\hat u_{t,\epsilon}$ to all of $M$ as
\[
u_{t,\epsilon}:=\begin{cases}
\hat u_{t,\epsilon} & \textrm{in }B_\epsilon,\\
u_t & \textrm{in }M\setminus B_\epsilon,
\end{cases}
\]
it turns out that 
\[
E_p(u_{t,\epsilon})<E_p(u_t). 
\]
Moreover, for $\epsilon$ small enough $u_{t,\epsilon}\in \mathcal H_{u_t}=\mathcal H_u$ and, as it is clear from the proof of Theorem \ref{th_Wei-existence}, it must be 
\[
E_p(u_{t,\infty})\leq E_p(u_{t,\epsilon})<E_p(u_t)=E_p(u).
\]
This contradicts \eqref{energies2} and concludes the proof.
\end{proof}

\begin{remark}\label{compact}\rm{
The $p$-harmonic general comparison theorem we have just proved does not recover completely Theorem \ref{Ith_SY_Topo}. First, due to the bad behaviour of the $p$-tension field with respect to the Riemannian product, here the map $M\times\mathbb R\to N$ given by $(q,t)\mapsto u_t(q)$ in general is not $p$-harmonic.\\
Moreover, in Theorems \ref{Ith_SY_Topo} the maps $u_t$ giving the homotopy are proven to be harmonic also when the target manifold $N$ is non-compact. This is achieved taking advantage of the solution to the Dirichlet problem for maps from bounded domains in $M$ to non-compact non-positively curved targets $N$. Indeed, developing the heat flow method used by R. Hamilton, \cite{Ham}, this result was achieved by Schoen and Yau in Section IX.8 of \cite{SY-LN}. To the best of our knowledge, in the $p$-harmonic setting, the Dirichlet problem for maps to non-compact manifolds has not been faced yet.   }
\end{remark}


\begin{thebibliography}{99}



%
%
%
%
\bibitem[A]{A} T. Aubin, \textit{Some nonlinear problems in Riemannian geometry.} Springer Monographs in Mathematics. Springer-Verlag, Berlin, 1998. xviii+395 pp.
%
%

\bibitem[B] {Burstall-London}F. Burstall, \textit{Harmonic maps of finite energy
from non-compact manifolds.} Jour. London Math. Soc. \textbf{30} (1984), 361--370.
%
%
%
%
%

\bibitem[CE]{CE} J. Cheeger, D.G. Ebin, \textit{Comparison theorems in Riemannian geometry.} North-Holland Mathematical Library, Vol. \textbf{9}. North-Holland Publishing Co., Amsterdam-Oxford; 1975. 
%
%
%
%
%
\bibitem[EL]{EL} J. Eells, L. Lemaire, \textit{Selected topics in harmonic maps.} CBMS Regional Conference Series in Mathematics, \textbf{50} (1983).

\bibitem[ES]{ES} J. Eells, J.H. Sampson, \textit{Harmonic mappings of Riemannian manifolds.}  Amer. J. Math.  \textbf{86}  1964 109--160.
%
%

%
\bibitem[GT] {GT}V. Gol'dshtein, M. Troyanov,\textit{\ The
Kelvin-Nevanlinna-Royden criterion for }$p$\textit{-parabolicity.} Math Z.
\textbf{232} (1999), 607--619.
%
%
%

%
\bibitem[Ham]{Ham} R. Hamilton, \textit{Harmonic maps of manifolds with boundary.} Lecture notes, Mathematics, No. \textbf{471} Springer, Berlin, Heidelberg, New York (1975). 

\bibitem[HL]{HL} R. Hardt, F.-H. Lin, \textit{Mappings minimizing the }$L^p$\textit{ norm of the gradient.}  Comm. Pure Appl. Math.  \textbf{40}  (1987),  no. 5, 555--588.

\bibitem[Har]{Ha} P. Hartman, \textit{On homotopic harmonic maps.} Canad. J. Math. \textbf{19} 1967 673--687.
%

\bibitem[Ho]{Holo-Dissertation} I. Holopainen, \textit{Nonlinear potential theory and quasiregular mappings on Riemannian manifolds.}
Ann. Acad. Sci. Fenn. Ser. A I Math. Dissertationes,  \textbf{74} (1990), 45 pp.
%
%
%

\bibitem[HPV]{HPV} I. Holopainen, S. Pigola, G. Veronelli, \textit{Global comparison principles for the $p$-Laplace operator on Riemannian manifolds.} To appear on Potential Analysis.
%
%

\bibitem [K]{Kawai}S. Kawai, $\mathit{p}$\textit{-harmonic maps and convex
functions.} Geom. Dedicata \textbf{74} (1999), 261--265.
%
%
%
%
%

\bibitem [LS]{LS}T. Lyons, D. Sullivan,\textit{\ Function theory, random paths and
covering spaces.} J. Diff. Geom. \textbf{18} (1984), 229--323.
%
%
%
%
%
%
%

%

\bibitem [PRS]{PRS-MathZ}S. Pigola, M. Rigoli, A.G. Setti, \textit{Constancy of
p-harmonic maps of finite q-energy into non-positively curved manifolds.} Math. Z. \textbf{258} (2008), no. 2, 347--362.
%
%

\bibitem[PST]{PST}S. Pigola, A.G. Setti, M. Troyanov. \textit{The topology at infinity of a manifold supporting an }$L^{p,q}$\textit{-Sobolev inequality.} Preprint.
%
\bibitem[PV]{PV}S. Pigola, G. Veronelli, \textsl{On the homotopy class of maps with finite $p$-energy into non-positively curved manifolds.} Geom. Dedic. \textbf{143} (2009), 109–-116. 
%
%
\bibitem[SU]{SU} R. Schoen, K. Uhlenbeck, \textit{A regularity theory for harmonic maps.}  J. Differential Geom. \textbf{17}  (1982), no. 2, 307--335.

\bibitem [SY1]{SY-Helv}R. Schoen and S.T. Yau, \textit{ Harmonic Maps and the
Topology of Stable Hypersurfaces and Manifolds with Non-negative Ricci
Curvature}. Comm. Math. Helv. \textbf{51} (1976), 333-341.

\bibitem [SY2]{SY-Topo}R. Schoen, S.T. Yau, \textit{Compact group actions and the
topology of manifolds with nonpositive curvature.} Topology \textbf{18}
(1979), 361--380.

\bibitem [SY3]{SY-LN}R. Schoen, S.T. Yau,\textit{\ Lectures on Harmonic Maps}.
Lectures Notes in Geometry and Topology, Volume II, International Press, Cambridge, MA, 1997. vi+394 pp.

%
\bibitem[S]{Sp}E.H. Spanier, \textit{Algebraic topology.} Springer-Verlag, New York, 1966. xvi+528 pp.
%
%

\bibitem[T]{Troy}M. Troyanov, \textsl{Parabolicity of manifolds. }
Siberian Adv. Math. \textbf{9} (1999), no. 4, 125--150.


\bibitem[VV]{VV} D. Valtorta, G. Veronelli, \textit{Stokes' theorem, volume growth and parabolicity.} Submitted.

\bibitem[V]{V} G. Veronelli, \textit{On $p$-harmonic maps and convex functions.} Manuscripta Math. \textbf{131} (2010), no. 3-4, 537--546. 

\bibitem [We1]{Wei-conference}S.W. Wei, \textit{The minima of the p-energy
functional. }Elliptic and parabolic methods in geometry (Minneapolis 1994),
pp. 171--203. A K Peters, Wellesley (1996).

\bibitem [We2]{Wei-indiana}S.W. Wei, \textit{Representing homotopy groups and
spaces of maps by }$\mathit{p}$\textit{-harmonic maps. }Indiana Univ. Math. J.
\textbf{47} (1998), 625--670.
%
%
%
%
\bibitem[Wh]{Wh} B. White, \textit{Homotopy classes in Sobolev spaces and the existence of energy minimizing maps.} Acta Math. \textbf{160} (1988), no. 1-2, 1--17.
%
\end{thebibliography}
\end{document}